\documentclass[reqno, 11pt]{amsart}
\usepackage{amssymb,graphicx}
\usepackage[usenames,dvipsnames]{color}
\usepackage{ifpdf}
\usepackage[pdftex]{thumbpdf}       

\ifpdf
\usepackage[pdftex,colorlinks]{hyperref} 
\hypersetup{%
pdftitle={A Lindley-type equation arising from a carousel problem},
pdfauthor={M. Vlasiou, I.J.B.F. Adan, J. Wessels},
pdfkeywords={},
bookmarksnumbered,
pdfstartview={FitH},
urlcolor=cyan,
}%
\else
  \newcommand\phantomsection\relax
  \newcommand{\url}[1]{#1}
  \newcommand{\href}[2]{#2}
\fi

\newcommand{\m}[1]{\mathcal{#1}}
\newcommand{\suml}{\sum\limits}
\newcommand{\intl}{\int\limits}
\newcommand{\e}{\mathbb{E}}
\newcommand{\p}{\mathbb{P}}

\theoremstyle{plain}
\newtheorem{theorem}{Theorem}

\newtheorem{corollary}{Corollary}

\newtheorem{lemma}{Lemma}

\theoremstyle{remark}
\newtheorem*{remark}{Remark}

\numberwithin{equation}{section}

\begin{document}
\title[A Lindley-type equation]{A Lindley-type equation\\ arising from a carousel problem}
\author[M. Vlasiou {\em et al.}]{M. Vlasiou}
\author[]{J. Wessels}
\author[]{I.J.B.F. Adan}

\address[M. Vlasiou, I.J.B.F. Adan and J. Wessels]{EURANDOM\\ P.O. Box 513\\ 5600 MB Eindhoven\\ The Netherlands.}
\address[I.J.B.F. Adan]{Department of Mathematics \& Computer Science\\ Eindhoven University of Technology\\ P.O. Box 513\\ 5600 MB Eindhoven\\ The Netherlands.}

\email[M. Vlasiou]{vlasiou@eurandom.tue.nl}
\email[I.J.B.F. Adan]{iadan@win.tue.nl}

\date{November 7, 2003}

\begin{abstract}
In this paper we consider a system with two carousels operated by one picker. The items to be picked are randomly located on the carousels and the pick times follow a phase-type distribution. The picker alternates between the two carousels, picking one item at a time. Important performance characteristics are the waiting time of the picker and the throughput of the two carousels. The waiting time of the picker satisfies an equation very similar to Lindley's equation for the waiting time in the $PH/U/1$ queue. Although the latter equation has no simple solution, we show that the one for the waiting time of the picker can be solved explicitly. Furthermore, it is well known that the mean waiting time in the $PH/U/1$ queue depends on to the complete interarrival time distribution, but numerical results show that, for the carousel system, the mean waiting time and throughput are rather insensitive to the pick-time distribution.

\end{abstract}

\maketitle

\section{Introduction}

In this paper we shall explore various methods to analyse a Lindley-type equation that emerges from a model involving two carousels alternately served by a picker. This equation differs from the original Lindley equation only in the change of a plus sign into a minus sign. The implications of this minor difference are rather far reaching, since in our situation there is an explicit solution and the result is surprisingly simple, while Lindley's equation has no simple solution. Furthermore, numerical results show that in this carousel model the mean waiting time is not very sensitive to the coefficient of variation of the pick time, which is in complete contrast to Lindley's equation.

Before getting into the details of the model, we describe the basic characteristics of carousels. A carousel is an automated storage and retrieval system, widely used in modern warehouses. It consists of a number of shelves or drawers rotating in a closed loop and it is operated by a picker that has a fixed position in front of the carousel. Carousels come in a huge variety of configurations, sizes and types. They can be horizontal or vertical and rotate in either one or both directions. Carousels are used in many different situations. For example,  e-commerce companies use them to store small items and manage small individual orders.

Carousel models have received much attention in the literature and continue to pose interesting problems. Jacobs {\it et al.} \cite{jacobs00}, for example, assumed a fixed number of orders and proposed a heuristic defining how many pieces of each item should be stored on the carousel in order to maximise the number of orders that can be retrieved without reloading. Usually a carousel is modelled as a circle. Stern \cite{stern86} and Ghosh and Wells \cite{ghosh92} considered a discrete model, where the circle consists of a fixed number of locations. Bartholdi and Platzman \cite{bartholdi86} and van den Berg \cite{vdberg96} proposed a continuous version, where the circle has unit length and the locations of the required items are represented as arbitrary points on the circle. In \cite{bartholdi86} the authors were mainly concerned with sequencing batches of requests in a bidirectional carousel, while in \cite{vdberg96} multiple order pick sequencing was studied. Ha and Hwang \cite{ha94} showed that performance is improved when some assignments of items to the set of drawers are more likely than others. Rouwenhorst {\it et al.}\ \cite{rouwenhorst96} gave stochastic upper bounds for the minimum travel time and studied the distribution of the travel time, assuming that the carousel changes direction after collecting at most one item. Litvak and Adan \cite{litvak01} and Litvak {\it et al.}\ \cite{litvak01a} assumed that the positions of the items are independent and uniformly distributed and gave a detailed analysis of the nearest-item heuristic, in which the next item to be picked is always the nearest one. More recent literature includes the work of Wan and Wolff \cite{wan04} that focused on minimising the travel time for `clumpy' orders and introduced the nearest-endpoint heuristic for which they obtained conditions for it to be optimal.

While almost all work concerns one-carousel models, real applications have triggered the study of models involving more complicated systems. Emerson and Schmatz \cite{emerson81} studied different storage schemes in a two-carousel setting by using simulation models. Recently, Hassini and Vickson \cite{hassini03} studied storage locations for items to minimise the long-run expected travel time in a two-carousel setting, while Park {\it et al.}\ ~\cite{park03} derived, under specific assumptions for the pick times, the distribution of the waiting time of the picker that alternates between two carousels. This allowed them to derive expressions for the system throughput and the picker utilisation. Our paper is motivated by their work. We extend the model by allowing for more general distributions for the pick times than those studied in \cite{park03} and we propose a different approach to the problem, leading to more explicit results.

In Section \ref{s:model} we will introduce the model in detail and analyse various implications of the nonstandard sign in the Lindley-type equation. In Section \ref{s:erlang} we will first consider pick times with an Erlang distribution and prove that the density of the waiting time of the picker can be expressed as a sum of exponentials. In Section \ref{s:PH} we extend this result to pick times with a phase-type distribution. In Section \ref{s:NumericalResults} we discuss some numerical results demonstrating that the throughput is fairly insensitive to the squared coefficient of variation of the pick times; the dominant factor is just the mean. We conclude with a brief summary and further research plans in Section \ref{s:conclusion}.

\section{The model} \label{s:model}

We consider a system consisting of two identical carousels and one picker. At each carousel there is an infinite supply of pick orders that need to be processed. The picker alternates between the two carousels, picking one order at a time. An important performance characteristic  is the throughput, i.e.\ the number of orders processed per unit time. Park {\em et al.}\ ~\cite{park03} determined the throughput when the pick times are either deterministic or exponentially distributed. We consider pick times following a phase-type distribution and derive explicit expressions for the throughput. Phase-type distributions may be used to approximate any pick-time distribution; see Schassberger~\cite{schassberger-W}.

Following Park {\em et al.}\ ~\cite{park03} we model a carousel as a circle of length $1$ and assume that it rotates in one direction at unit speed. Each pick order requires exactly one item. The picking process may be visualised as follows. When the picker is about to pick an item at one of the carousels, he may have to wait until the item is rotated in front of him. In the meantime, the other carousel rotates towards the position of the next item. After completion of the first pick the carousel is instantaneously replenished and the picker turns to the other carousel, where he may have to wait again, and so on. Let the random variables $P_n$, $R_n$ and $W_n$ ($n \geqslant 1$) denote the pick time, rotation time and waiting time for the $n$th item. Clearly, the waiting times $W_n$ satisfy the recursion
\begin{equation}\label{lindley1}
W_{n+1}=(R_{n+1}-P_n-W_n)^+, \qquad n = 0, 1, \ldots; \qquad P_0 = W_0 \stackrel{\mathsf{def}}{=}0,
\end{equation}
where $(x)^+=\max\{0, x\}.$ We assume that both $\{P_n, n \geqslant 1\}$ and $\{R_n, n \geqslant 1\}$ are sequences of independent identically distributed  random variables, also independent of each other. The pick times $P_n$ have a phase-type distribution $G(\cdot)$ and the rotation times $R_n$ are uniformly distributed on $[0,1)$ (which means that the items are randomly located on the carousels). Then $\{W_n\}$ is a Markov chain, with state space $[0,1)$.
In \cite{park03} it is shown that $\{W_n\}$ is an aperiodic, recurrent Harris chain, which possesses a unique equilibrium distribution.
In equilibrium, equation (\ref{lindley1}) becomes
\begin{equation} \label{recursion}
W\stackrel{\m{D}}{=}(R-P-W)^+.
\end{equation}
Note the striking similarity to Lindley's equation for the waiting times in a single-server queue: the only difference is the sign of $W_n$. Let $\pi_0 =\p[W=0]$ and $f(\cdot)$ denote the density of $W$ on $[0,1]$. From \eqref{recursion} it readily follows that (cf. Equation (3) in Park {\em et al.}\ ~\cite{park03})
\begin{equation} \label{e3}
f(x)=\pi_0G(1-x)+\int_0^{1-x}G(1-x-z)f(z)dz, \qquad 0 \leqslant x \leqslant 1,
\end{equation}
with the normalisation equation
\begin{equation}
\label{norm}
\pi_0 + \int_0^1 f(x) dx = 1.
\end{equation}
Once the solution to equations \eqref{e3} and \eqref{norm} is known, we can compute $\e[W]$ and thus also the throughput $\tau$ from
\begin{equation}
\label{tau}
\tau = \frac{1}{\e[W]+\e[P]}.
\end{equation}

As pointed out before, \eqref{recursion} (with a plus sign instead of minus sign for $W$) is precisely Lindley's equation for the stationary
waiting time in a $PH/U/1$ queue. This equation has no simple solution, but we show that the waiting time of the picker can be solved for explicitly. Lindley's equation is one of the most studied equations in queueing theory. For excellent textbook treatments we refer to Asmussen \cite{asmussen-APQ}, Cohen \cite{cohen-SSQ}, and the references therein. It is interesting to investigate the impact on the analysis of such a slight modification to the original equation.\\

In the following we explore various methods of solving the Lindley-type recursion \eqref{recursion}, or equivalently \eqref{e3} and \eqref{norm}. Since  \eqref{e3} is a Fredholm type equation, a natural way to proceed is by successive substitutions. This yields the formal solution
\begin{equation} \label{formal}
f(x)=\pi_0 \sum_{j=1}^{\infty} G^{j\,*}(1-x), \qquad 0 \leqslant x \leqslant 1,
\end{equation}
where
\[
G^{1\,*}(1-x) \stackrel{\mathsf{def}}{=} G(1-x); \qquad G^{n\,*}(1-x)\stackrel{\mathsf{def}}{=} \int_0^{1-x}G(1-x-z)\,G^{(n-1)\,*}(1-z)\,dz, \qquad n \geqslant 2.
\]
Since $G(\cdot)$ is a distribution, from the last relation we have, for $n\geqslant 1,$ that
$$
G^{(n+2)\,*}(x)\leqslant\int_0^x G^{(n+1)\,*}(1-z)\,dz\leqslant\int_0^x\int_0^{1-z}G^{n\,*}(1-y)\,dydz=\int_0^x\int_z^1G^{n\,*}(y)\,dydz,
$$
which implies that $G^{3\,*}(x)\leqslant 1/2.$ Now, by induction, it can be easily shown that, for $n\geqslant 1$
\[
G^{2(n+1)\,*}(x) \leqslant G^{(2n+1)\,*}(x) \leqslant \frac{1}{2^n} , \qquad 0 \leqslant x \leqslant 1.
\]
This means that the infinite sum \eqref{formal} converges (uniformly) for $0 \leqslant x \leqslant 1.$

However, for a non-trivial distribution $G(\cdot)$, one cannot easily compute $f(\cdot)$ using \eqref{formal}. The difficulty lies in the fact that $G^{n\, *}(\cdot)$ is not the $n$-fold convolution of the distribution function $G(\cdot)$. Therefore, we need a method that leads to more tractable results. For this reason we proceed by applying Laplace transforms to solve \eqref{e3}. Laplace transforms are a standard approach for solving the original Lindley equation. For our model, this approach yields explicit and computable expressions for the density $f(\cdot)$ and the throughput $\tau$, involving roots of a certain equation.

Another possibility is to obtain from \eqref{e3} a solvable differential equation. This method was used to some extent in Park {\it et al.} \cite{park03}. They focused on deterministic and exponentially distributed pick times and commented that ``the approach of deriving a differential equation for each pick-time distribution was rather ad hoc". However, this method can be generalised to include phase-type distributions as well. For more details we refer to Vlasiou {\it et al.}\ \cite{vlasiou03}. The advantage of the Laplace transform approach over this one is that it leads to a more explicit solution.

\section{Erlang pick times} \label{s:erlang}
In this section we will use Laplace transforms to solve \eqref{recursion} and we compare this method to the work that has previously been done in \cite{park03}. Throughout this section we assume that the pick times follow an Erlang distribution Erl$(\mu, n)$ with scale parameter $\mu$ and $n$ stages, that is
\[
G(x)=1-e^{-\mu x}\sum_{j=0}^{n-1}\frac{(\mu x)^j}{j!}, \qquad x \geqslant 0.
\]
Let $\phi (\cdot)$ denote the Laplace transform of $f(\cdot)$ over the interval $[0, 1]$, i.e.\
\[
\phi(s)=\intl_0^1e^{-sx}f(x)dx.
\]
We emphasise that, for the Laplace transform over a bounded interval, the standard properties are no longer valid, in the sense that there are no standard results for calculating the inverse transform over a bounded interval. Note that $\phi(\cdot)$ is analytic in the whole complex plane. It is convenient to replace $x$ by $1-x$ in \eqref{e3}, yielding
\begin{equation}\label{1-x}
f(1-x)=\pi_0G(x)+\int_0^{x}G(x-z)f(z)dz, \qquad 0 \leqslant x \leqslant 1.
\end{equation}
By taking the Laplace transform of \eqref{1-x} and using \eqref{norm} we obtain
\begin{align*}
e^{-s}\phi(-s)&=\pi_0\left(\frac{1-e^{-s}}{s}-\sum_{j=0}^{n-1} \frac{\mu^j}{(\mu+s)^{j+1}}+ \sum_{j=0}^{n-1}\sum_{i=0}^j \frac{\mu^j}{i!(\mu+s)^{j+1-i}}e^{-(\mu+s)}\right)\\
&\qquad -\frac{e^{-s}}{s}(1-\pi_0)+\frac{1}{s}\phi(s)-\sum_{j=0}^{n-1} \frac{\mu^j}{(\mu+s)^{j+1}}\phi(s)\\
&\qquad +e^{-(\mu+s)}\sum_{j=0}^{n-1}\sum_{i=0}^j\sum_{\ell=0}^i \binom{i}{\ell}\frac{\mu^j}{i!(\mu+s)^{j+1-i}}\phi^{(\ell)}(-\mu),
\end{align*}
which, by rearranging terms and using the identity
\[
\sum_{j=0}^{n-1} \frac{\mu^j}{(\mu+s)^{j+1}} =
\frac{(\mu+s)^n - \mu^n}{s (\mu+s)^n} ,
\]
can be simplified to
\begin{align}\label{erl_iter}\nonumber
e^{-s}\phi(-s)
- \frac{\mu^{n}}{s(\mu+s)^{n}} \, \phi(s)
&=\pi_0 \left(\frac{\mu^{n}}{s(\mu+s)^{n}}+ e^{-(\mu+s)} \sum_{j=0}^{n-1}\sum_{i=0}^j \frac{\mu^j}{i!(\mu+s)^{j+1-i} }\right)\\
&\qquad-\frac{e^{-s}}{s}+e^{-(\mu+s)}\sum_{j=0}^{n-1}\sum_{i=0}^j\sum_{\ell=0}^i \binom{i}{\ell} \frac{\mu^j}{i!(\mu+s)^{j+1-i}}\phi^{(\ell)}(-\mu).
\end{align}
In the above expression, $\phi^{(\ell)}(\cdot)$ denotes the $\ell$th derivative of $\phi(\cdot)$. Note that both $\phi(-s)$ and $\phi(s)$ appear in \eqref{erl_iter}. To obtain an additional equation we replace $s$ by $-s$ in \eqref{erl_iter} and form a system from which $\phi(s)$ can be solved, yielding the following theorem.

\begin{theorem}
For all $s$, the transform $\phi(s)$ satisfies
\begin{equation}\label{erl_polywn}
\phi(s) R(s) = -e^{-s}s(\mu+s)^{n}A(-s)-\mu^{n}A(s),
\end{equation}
where
\begin{eqnarray*}
R(s)&=&
s^2(\mu^2-s^2)^{n}+\mu^{2n},\\
A(s)&=&\pi_0\left(\mu^{n}+e^{-(\mu+s)}\sum_{j=0}^{n-1}\sum_{i=0}^j\frac{s\mu^j(\mu+s)^{n-j-1+i}}{i!}\right)-e^{-s}(\mu+s)^{n}\\
&&\qquad+e^{-(\mu+s)}\sum_{j=0}^{n-1}\sum_{i=0}^j\sum_{\ell=0}^i\binom{i}{\ell}\frac{s\mu^j(\mu+s)^{n-j-1+i}}{i!}\phi^{(\ell)}(-\mu).
\end{eqnarray*}
\end{theorem}

In \eqref{erl_polywn} we still need to determine the $n+1$ unknowns $\pi_0$ and $\phi^{(\ell)}(-\mu)$ for $\ell=0,\ldots,n-1$. Note that for any zero of the polynomial $R(\cdot)$, the left-hand side of \eqref{erl_polywn} vanishes (since $\phi(\cdot)$ is analytic everywhere). This implies that the right-hand side should also vanish. Hence, the zeros of $R(\cdot)$ provide the equations necessary to determine the unknowns.

\begin{lemma}\label{lem:zero}
The polynomial $R(\cdot)$ has exactly $2n+2$ simple zeros $r_1, \ldots, r_{2n+2}$ satisfying $r_{2n+3-i} = -r_{i},$ for $i = 1, \ldots, n+1$.
\end{lemma}
\begin{proof} Since $R(s)$ is a polynomial in $s^2$ of degree $n+1,$ it follows that $R(s)$ has exactly $2n+2$ zeros, with the property that
each zero $s$ has a companion zero $-s$. Furthermore, it is easily verified that $\gcd[R(s),R^{\prime}(s)]=1$. This means that the polynomials $R(s)$ and $R^{\prime}(s)$ have no common factor of degree greater than zero, or that $R(s)$ has only simple zeros. \end{proof}

In the following lemma we prove that the $2n+2$ zeros of $R(\cdot)$ produce $n+1$ independent linear equations for the unknowns.

\begin{lemma}\label{lem:eqs}
The probability $\pi_0$ and the quantities $\phi^{(\ell)}(-\mu),$  $\ell=0,\ldots,n-1$ are the unique solution to the $n+1$ linear equations,
\begin{equation*}
e^{-r_i}r_i(\mu+r_i)^{n}A(-r_i)+\mu^{n}A(r_i) = 0, \qquad i = 1, \ldots, n+1.
\end{equation*}
\end{lemma}
\begin{proof} For any zero of $R(\cdot)$ the right-hand side of \eqref{erl_polywn} should vanish. Hence, for two companion zeros $r_i$ and $r_{2n+3-i} = -r_i$, $i = 1, \ldots, n+1$, we have
\begin{eqnarray}
\label{a1}
e^{-r_i}r_i(\mu+r_i)^{n}A(-r_i)+\mu^{n}A(r_i)  & = & 0,\\
\label{a2}
-e^{r_i}r_i(\mu-r_i)^{n}A(r_i)+\mu^{n}A(-r_i)  & = & 0.
\end{eqnarray}
The determinant of \eqref{a1} and \eqref{a2}, treated as equations for $A(-r_i)$ and $A(r_i)$, is equal to $R(r_i) = 0$. Hence, \eqref{a1} and \eqref{a2} are dependent, and so we may omit one of them. This leaves a system of $n+1$ linear equations for the unknowns $\pi_0$ and
$\phi^{(\ell)}(-\mu),$ $\ell=0,\ldots,n-1$.
The uniqueness of the solution follows from the general theory of Markov chains that implies that there is a unique equilibrium distribution and thus also a unique solution to \eqref{erl_iter}. \end{proof}

Once $\pi_0$ and $\phi^{(\ell)}(-\mu)$, $\ell=0,\ldots,n-1$ are determined, the transform $\phi(\cdot)$ is known. It remains to invert the transform. By collecting the terms that include $e^{-s}$ we can rewrite \eqref{erl_polywn} in the form
\begin{equation} \label{erl_decompose}
\phi(s) = \frac{P(s)}{R(s)}+e^{-s}\frac{Q(s)}{R(s)},
\end{equation}
where $P(s)$ and $Q(s)$ are polynomials of degree $2n+1$ and $n+1$ respectively. Note that, without the last term, the transform is  rational so the inverse would be straightforward if we had Laplace transforms on $[0, \infty)$. As it is, we must proceed more carefully. Since $\deg[R]$ is greater than $\deg[P]$ and $\deg[Q]$, \eqref{erl_decompose} can be decomposed into distinct irreducible fractions. This leads to
\[
\phi(s)= \frac{c_1}{s-r_1}+\cdots+\frac{c_{2n+2}}{s-r_{2n+2}}+e^{-s} \left[\frac{\hat{c}_1}{s-r_1} +\cdots+\frac{\hat{c}_{2n+2}}{s-r_{2n+2}}\right],
\]
where the coefficients $c_i$ and $\hat{c}_i$ are given by
\begin{equation} \label{ci}
c_i = \lim_{s\to r_i}\frac{P(s)}{R(s)} \, (s-r_i) = \frac{P(r_i)}{R^\prime(r_i)}, \qquad \hat{c}_i = \lim_{s\to r_i}\frac{Q(s)}{R(s)} \, (s-r_i)=
\frac{Q(r_i)}{R^\prime(r_i)}.
\end{equation}
Note that the derivative $R^\prime(r_i)$ is nonzero, since $r_i$ is a simple zero. Since $\phi(s)$ is analytic everywhere, then for every root $r_i$ of $R(s)$ we have
\[
P(r_i) = -e^{-r_i} Q(r_i) , \qquad i = 1, \ldots, 2n+2.
\]
Hence, from \eqref{ci} it follows that
\begin{equation}\label{relation}
c_i=-e^{-r_i}\hat{c}_i,
\end{equation}
and thus
\[
\phi(s)=\sum_{i=1}^{2n+2}\frac{c_i}{s-r_i}\left[1-e^{r_i-s}\right],
\]
which is the transform (over a bounded interval) of a mixture of $2n+2$ exponentials. Now that the density is known, \eqref{norm} can be used to derive a simple explicit expression for $\pi_0$. These findings are summarised in the following theorem.

\begin{theorem} \label{th:density}
The density of $W$ on $[0,1]$ is given by
\begin{equation}\label{density}
  f(x)=\suml_{i=1}^{2n+2}c_i e^{r_i x}, \qquad 0 \leqslant x \leqslant 1,
\end{equation}
and
\begin{equation}\label{pi0}
\pi_0 = \p[W=0] = 1+\sum_{i=1}^{2n+2} \frac{c_i}{r_i}(1-e^{r_i}).
\end{equation}
\end{theorem}

\begin{corollary}\label{cor:tau}
The throughput $\tau$ satisfies
\[
\tau^{-1} = \e[P] + \e[W] =
\frac{n}{\mu} +
\suml_{i=1}^{2n+2} \frac{c_i}{r^2_i} [1 + (r_i -1)e^{r_i}] .
\]
\end{corollary}

Although the roots $r_i$ and coefficients $c_i$ may be complex, the expressions \eqref{density} and \eqref{pi0} will be positive. This follows from the fact that the equilibrium equation \eqref{e3} and the normalisation equation \eqref{norm} have a unique solution. Of course, it is also clear that each root $r_i$ and coefficient $c_i$ have a companion conjugate root and conjugate coefficient, which implies that the imaginary parts in \eqref{density} and \eqref{pi0} cancel.

\section{Phase-Type pick times}\label{s:PH}
Let us now assume that the pick times follow an Erl$(\mu, n)$ with probability $\alpha_n$, $n = 1, \ldots, N.$ In other words,
\begin{equation} \label{mixed}
G(x)=\sum_{n=1}^N \alpha_n\left(1-e^{-\mu x}\sum_{j=0}^{n-1}\frac{(\mu x)^j}{j!}\right) ,\qquad x \geqslant 0.
\end{equation}
The class of the phase-type distributions of the above form is dense in the space of distribution functions defined on $[0, \infty)$. This means that for any such distribution function $F(\cdot)$, there is a sequence $F_n (\cdot)$ of phase-type distributions of this class  that converges weakly to $F(\cdot)$ as $n$ goes to infinity; for details see Schassberger~\cite{schassberger-W}. Below we give the result for pick time distributions of the form \eqref{mixed}.\\

The analysis proceeds along the same lines as in Section \ref{s:erlang}. The formulae in the intermediate steps are simply linear combinations of the ones that appear for Erlang pick times. This leads to the following result.
\begin{theorem}
For all $s$, the transform $\phi(s)$ satisfies
\begin{equation}\label{ph_polywn}
\phi(s) \widetilde{R}(s) =
-e^{-s}s(\mu+s)^{N}\widetilde{A}(-s)-\sum_{n=1}^N\alpha_{n}\mu^{n}(\mu-s)^{N-n}\widetilde{A}(s),
\end{equation}
where
\begin{eqnarray*}
\widetilde{R}(s) & = &
s^2(\mu^2-s^2)^{N}+\sum_{n=1}^N\sum_{m=1}^N\alpha_{n}\alpha_{m}\mu^{n}\mu^{m}(\mu-s)^{N-n}(\mu+s)^{N-m},\\
\widetilde{A}(s) & = &
\pi_0\sum_{n=1}^N\alpha_{n}\left(\mu^{n}(\mu+s)^{N-n}+e^{-(\mu+s)}\sum_{j=0}^{n-1}\sum_{i=0}^j\frac{s\mu^j(\mu+s)^{N-j-1+i}}{i!} \right)\\
&& \qquad +\sum_{n=1}^N\alpha_{n}\left(-e^{-s}(\mu+s)^{N}+e^{-(\mu+s)}\sum_{j=0}^{n-1}\sum_{i=0}^j\sum_{\ell=0}^i\binom{i}{\ell}\frac{s\mu^j(\mu+s)^{N-j-1+i}}{i!} \phi^{(\ell)}(-\mu)\right).
\end{eqnarray*}
\end{theorem}

The unknowns $\pi_0$ and $\phi^{(\ell)}(-\mu),$ $\ell=0,\ldots,n-1$ can be determined in the same way as in Section \ref{s:erlang}. The polynomial $\widetilde{R}(\cdot)$ has exactly $2N+2$ zeros, with the property that each zero $s$ has a companion zero $-s$. We {\em assume} that all these zeros are simple and label them $\widetilde{r_1}, \ldots, \widetilde{r}_{2N+2}$ such that $\widetilde{r}_{2N+3-i} = - \widetilde{r}_i$ for $i = 1,\ldots, N+1$.
Then the following lemma can be readily established.

\begin{lemma}\label{lem:ph_eqs}
The probability $\pi_0$ and the quantities $\phi^{(\ell)}(-\mu),$ $\ell=0,\ldots,n-1$ are the unique solution to the $N+1$ linear equations,
\begin{equation}\label{ph_eqA}
e^{-\widetilde{r}_i}\widetilde{r}_i(\mu+\widetilde{r}_i)^{N}\widetilde{A}(-\widetilde{r}_i)+\sum_{n=1}^N\alpha_{n}\mu^{n}(\mu-\widetilde{r}_i)^{N-n}
\widetilde{A}(\widetilde{r}_i) = 0, \qquad i = 1, \ldots, N+1.
\end{equation}
\end{lemma}

Given $\pi_0$ and $\phi^{(\ell)}(-\mu),$ $\ell=0,\ldots,n-1,$ the transform $\phi(\cdot)$ is completely known. Partial fraction decomposition of the transform yields
\[
\phi(s)=\sum_{i=1}^{2N+2}\frac{\widetilde{c}_i}{s-\widetilde{r}_i}\left[1-e^{\widetilde{r}_i-s}\right],
\]
from which we conclude that the density of the waiting time is a mixture of $2N+2$ exponentials. Hence, as was the case for Erlang pick times, the density is given by
\[
f(x)=\suml_{i=1}^{2N+2}\widetilde{c}_ie^{\widetilde{r}_ix}.
\]

\begin{remark}
When $R(\cdot)$ has multiple zeros, the analysis proceeds in essentially the same way. For example, if  $\widetilde{r}_1 = \widetilde{r}_2$ (so $\widetilde{r}_1$ and, thus, $\widetilde{r}_{2N+2}$ are double zeros), then \eqref{ph_eqA} is identical for $i=1$ and $i=2$. Nonetheless, an additional equation can be obtained by requiring that the derivative of the right-hand side of \eqref{ph_polywn} should vanish at $s=r_1$. The partial-fraction decomposition of $\phi(\cdot)$ then becomes
\begin{eqnarray*}
\phi(s) & = &
\frac{\widetilde{c}_1}{(s-\widetilde{r}_1)^2}
\left[1-e^{\widetilde{r}_1-s}-(s-\widetilde{r}_1)e^{\widetilde{r}_1-s}\right]+
\sum_{i=2}^{2N+1}
\frac{\widetilde{c}_i}{s-\widetilde{r}_i}
\left[1-e^{\widetilde{r}_i-s}\right] \\
&& +
\frac{\widetilde{c}_{2N+2}}{(s-\widetilde{r}_{2N+2})^2}
\left[1-e^{\widetilde{r}_{2N+2}-s}-(s-\widetilde{r}_{2N+2})e^{\widetilde{r}_{2N+2}-s}\right],
\end{eqnarray*}
the inverse of which is given by
\[
f(x)  =
\widetilde{c}_1 x e^{\widetilde{r}_1x} +
\suml_{i=2}^{2N+1}\widetilde{c}_ie^{\widetilde{r}_ix} +
\widetilde{c}_{2N+2} x e^{\widetilde{r}_{2N+2}x}.
\]
\end{remark}

\section{Numerical results}\label{s:NumericalResults}

\begin{figure}[htbp]
 \leavevmode
\begin{center}
\includegraphics{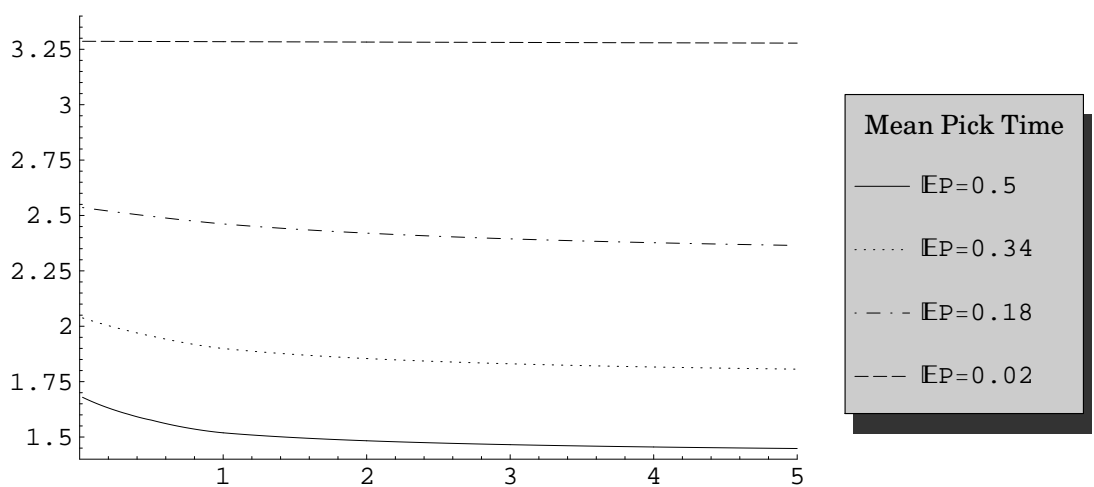}
\end{center}
\caption{Plot of throughput vs. the squared coefficient of variation
of the pick time.}
\label{fig1}
\end{figure}

This section is devoted to some numerical results. For various values of the mean pick time $\e[P]$ we plot in Figure \ref{fig1} the throughput $\tau$ versus the squared coefficient of variation of the pick time, $c^2_P$. The mean pick time is chosen to be comparable to the mean rotation time, which is $frac{1}{2}$. In each plot we fit a mixed Erlang or hyperexponential distribution to $\e[P]$ and $c^2_P$, depending on whether the squared coefficient of variation is less or greater than $1$ (see, for example, Tijms \cite{tijms-FCSM}).

Hyperexponential distributions form another useful class of phase-type distributions. They can be used to model pick times with squared coefficient of variation greater than 1. Furthermore, hyperexponential distributions are always unimodal, which is not the case for mixed Erlang distributions. The analysis for hyperexponential pick times is very similar to the one presented in the previous section.

So, if $1/n \leqslant c_P^2 \leqslant 1/(n-1)$ for some $n = 2, 3, \ldots$, then the mean and squared coefficient of variation of the mixed Erlang distribution
\[
G(x)= p
\left(1-e^{-\mu x}\sum_{j=0}^{n-2}\frac{(\mu x)^j}{j!}\right)
+ (1-p)
\left(1-e^{-\mu x}\sum_{j=0}^{n-1}\frac{(\mu x)^j}{j!}\right),
\qquad x \geqslant 0 ,
\]
matches with $\e[P]$ and $c_P^2$, provided the parameters $p$ and $\mu$ are chosen as
\[
p = \frac{1}{1+c_P^2}
[ n c_P^2 - \{ n(1+c_P^2) - n^2 c_P^2 \}^{1/2} ] , \qquad
\mu = \frac{n-p}{\e[P]} .
\]
On the other hand, if $c_P^2 > 1$, then the mean and squared coefficient of variation of the hyperexponential distribution
\[
G(x) = p_1 (1- e^{-\mu_1 x}) + p_2 (1-  e^{-\mu_2 x}) ,
\qquad x \geqslant 0,
\]
match with $\e[P]$ and $c_P^2$ provided the parameters $\mu_1, \mu_2, p_1$ and $p_2$ are chosen as
\begin{eqnarray*}
&& p_1 =
\frac{1}{2} \left( 1 + \sqrt{\frac{c_P^2 - 1}{c_P^2 + 1}}
\right) , \qquad p_2 = 1-p_1, \\
&& \mu_1 =
\frac{2 p_1}{\e[P]} \qquad \mbox{and}\qquad\mu_2 = \frac{2 p_2}{\e[P]}.
\end{eqnarray*}
For single-server queuing models it is well-known that the mean waiting time depends (approximately linearly) on the squared coefficients of variation of the interarrival (and service) times. The results in Figure \ref{fig1}, however, show that for the carousel model, the mean waiting time is not very sensitive to the squared coefficient of variation of the pick time and thus neither is the throughput $\tau$; it indeed decreases as $c^2_P$ increases, but very slowly.
This phenomenon may be explained by the fact that the waiting time of the picker is bounded by $1$, i.e.\ the time needed for a full rotation of the carousel.

\section{Concluding remarks and further research}\label{s:conclusion}

In this paper we have considered a system with two carousels operated by one picker. Using Laplace transforms over a bounded interval we have obtained an explicit solution for the density of the waiting time of the picker. We have shown that if we let the pick time follow a phase-type distribution, then the density is a mixture of exponentials. Numerical results show that the squared coefficient of variation of the pick time does not influence the throughput significantly.

We have solved the Lindley-type recursion \eqref{lindley1} under  specific assumptions on the random variables $R_n$ and $P_n$. In particular, we assumed that $R_n$ is uniformly distributed on $[0,1]$ and $P_n$ follows a phase-type distribution, for every $n$. This makes sense if one has a carousel application in mind. Nonetheless, it is mathematically interesting to try and solve this recursion under less restrictive assumptions. In further research we shall try to solve \eqref{lindley1} allowing $R_n$ and $P_n$ to follow a more general distribution.

\section*{Acknowledgements}
We would like to thank the referee for many helpful suggestions on the paper.

\bibliographystyle{abbrv}
\bibliography{maria}
\end{document}